\documentclass[12pt]{amsart}
\usepackage{amsmath,amssymb}
\usepackage[curve,matrix,arrow,cmtip]{xy}                       
\usepackage[english]{babel}
\usepackage[colorlinks,linkcolor=blue,anchorcolor=blue,citecolor=blue,backref=page]{hyperref}
\usepackage{tabulary}

\def\leq{\leqslant}
\def\geq{\geqslant}
\def\({\left(}
\def\){\right)}

\usepackage[norefs,nocites]{refcheck}
\long\def\beginFORGET#1\endFORGET{}

 \def\Gal{\mathop{\rm Gal}\nolimits}

\def\ord{\mathop{\rm ord}\nolimits}
\def\ind{\mathop{\rm ind}\nolimits}

\let\phi\varphi
\let\epsilon\varepsilon

\newtheorem{Thm}{Theorem}[section]
\newtheorem{Prop}[Thm]{Proposition}
\newtheorem{Lem}[Thm]{Lemma}

\newtheorem{Cor}[Thm]{Corollary}

\newtheorem{Rem}[Thm]{Remark}

\newtheorem*{theorem*}{Theorem}

\numberwithin{equation}{section}

\numberwithin{Thm}{section}

\def\qed{{\hskip0pt\unskip\unskip\nobreak\hfil\penalty50
          \hskip1em\hbox{}\nobreak\hfil
           {$\square$}
          \parfillskip=0pt\finalhyphendemerits=0
          \par}\medskip}

\newenvironment{Proof}
               {\noindent{\bf Proof.}\ }
               {\qed}

               {\noindent{\bf Proof of #1.}\ }
               {\qed}

\newcommand{\BC}{{\mathbb{C}}}

\newcommand{\BF}{{\mathbb{F}}}

\newcommand{\BQ}{{\mathbb{Q}}}

\newcommand{\BZ}{{\mathbb{Z}}}

\newcommand{\CO}{{\mathcal O}}

\newbox\mybox
\def\arrover#1{\mathrel{
       \setbox\mybox=\hbox spread 1.4em
              {\hfil$\scriptstyle#1$\hfil}
       \vbox{\offinterlineskip\copy\mybox
             \hbox to\wd\mybox{\rightarrowfill}}}}

\def\larrover#1{\mathrel{
       \setbox\mybox=\hbox spread 1.4em
              {\hfil$\scriptstyle#1\vphantom{g}$\hfil}
       \vbox{\offinterlineskip\copy\mybox
             \hbox to\wd\mybox{\leftarrowfill}}}}

\def\ontoover#1{\mathrel{
       \setbox\mybox=\hbox spread 1.4em
              {\hfil$\scriptstyle#1\vphantom{g}$\hfil}
       \vbox{\offinterlineskip\copy\mybox
             \hbox to\wd\mybox{\rightarrowfill\hskip-2.8mm
                               $\rightarrow$}}}}
\def\leftontoover#1{\mathrel{
       \setbox\mybox=\hbox spread 1.4em
              {\hfil$\scriptstyle#1\vphantom{g}$\hfil}
       \vbox{\offinterlineskip\copy\mybox
             \hbox to\wd\mybox{$\leftarrow$\hskip-2.8mm
                               \leftarrowfill}}}}

\let\onto\twoheadrightarrow

\begin{document}

\title[Multiplicative orders]{Multiplicative orders of Gauss periods and the arithmetic of real quadratic fields}
\author{Florian Breuer}
\address{School of Mathematical and Physical Sciences, University of Newcastle, Newcastle,NSW 2308,  Australia}
\email{florian.breuer@newcastle.edu.au}

\maketitle

\begin{abstract}
We obtain divisibility conditions on the multiplicative orders of elements of the form $\zeta + \zeta^{-1}$ in a finite field by exploiting a link to the arithmetic of real quadratic fields.
\end{abstract}

\section{Introduction}

Let $q$ be a prime number and $n$ a positive integer. We denote by $\BF_{q^n}$ the finite field of $q^n$ elements. 
Suppose $p=2n+1$ is an odd prime number and let $\zeta\in\BF_{q^{2n}}$ be a primitive $p^{\mathrm{th}}$ root of unity. We set
\[
\alpha = \zeta + \zeta^{-1}.
\]
Then $\alpha\in\BF_{q^n}$ is known as a Gauss period of type $(n,2)$ over $\BF_q$, and has many desirable properties. For example, when $\alpha$ is a primitive element, then it generates a normal basis for $\BF_{q^n}$. 
%This can only happen when $(\BZ/p\BZ)^* = \langle -1, q \rangle$, a situation we will assume in this paper. 
As a result, one is interested in the multiplicative order $\ord(\alpha)$ of $\alpha$ in $\BF_{q^n}^*$. 
See \cite{ASV10, GS98, GV95, Pop12, Pop14} and the references therein, where amongst other things lower bounds on $\ord(\alpha)$ are obtained.

In this paper, we will look at divisibility conditions, which imply upper bounds.
The trivial upper bound $\ord(\alpha)\leq q^n-1$ is often sharp when $q=2$ or $3$, and in general the index $\ind(\alpha) := (q^n-1)/\ord(\alpha)$ tends to be small. The goal of this paper is to show how certain small prime factors of this index can be detected in the arithmetic of the real quadratic field $\BQ(\sqrt{p})$.

 More precisely, denote by $\varepsilon_p$ and $h_p$ the fundamental unit and class number of $K=\BQ(\sqrt{p})$, respectively. Denote by $\CO_K$ the ring of integers of~$K$.
When $q$ is inert in $K/\BQ$, we denote by $\ind(\varepsilon_p \bmod q)$ the multiplicative index of $(\varepsilon_p \bmod q\CO_K)$ in $\left(\CO_K/q\CO_K\right)^* \cong \BF_{q^2}^*$.

%
%the residue class of $\varepsilon_p$ modulo the maximal ideal $q\CO_K$ in the finite field $\CO_K/q\CO_K$ of order $q^2-1$.

%
%
%$\varepsilon_p \bmod q \in (\CO_K/q\CO_K)^* \cong \BF_{q^2}^*$ the residue class of $\varepsilon_p$ modulo the maximal ideal $q\CO_K$, and by $\ind(\varepsilon_p \bmod q)$ its index in the cyclic group $(\CO_K/q\CO_K)^*$ of order $q^2-1$.

Our main result is the following.

\begin{Thm}\label{Thm:Main}
Let $p\equiv 5 \bmod 8$ be a prime number, suppose that $(\BZ/p\BZ)^* = \langle -1, q\rangle$
and let $\zeta\in\BF_{q^{p-1}}$ be a primitive $p^{\mathrm{th}}$ root of unity.
Then
\[
	\gcd\big(\ind(\zeta + \zeta^{-1}),q^2-1\big) = \ind(\varepsilon_p^{h_p} \bmod q).
\]
\end{Thm}

Related elements of interest are $\beta = \zeta + 1 \in \BF_{q^{2n}}$, whose multiplicative orders (when $q=2$) determine periods of Ducci sequences, see \cite{BLM07, Bre19, BS19}. When $q$ is a primitive root modulo $p$, 
then $p\nmid (q^n-1)$ so $p\nmid \ord(\zeta+\zeta^{-1})$ and we get
\[
\ord(\zeta+1) = \ord(\zeta^2+1) = \ord\big(\zeta(\zeta+\zeta^{-1})\big) = p\ord(\zeta + \zeta^{-1}),
\]
where we have used the fact that $\zeta$ and $\zeta^2$ are conjugate. 

Now, when $q=2$ and $p \equiv 1 \bmod 4$, we find that $(\BZ/p\BZ)^* = \langle -1, 2 \rangle$ is equivalent to $(\BZ/p\BZ)^* = \langle 2 \rangle$. Theorem \ref{Thm:Main} implies 

\begin{Cor}
Suppose that $p \equiv 5 \bmod 8$ is prime and that $2$ is a primitive root modulo $p$. Let $\zeta \in \BF_{2^{p-1}}$ be a primitive $p\mathrm{th}$ root of unity. 
Then the following are equivalent:
\begin{enumerate}
	\item $\ind(\zeta+1)$ is divisible by 3.
	\item $\ind(\zeta+\zeta^{-1})$ is divisible by 3.
	\item The eventual period $P$ of any Ducci sequence in $\BZ^p$ formed by iterating the map 
	\[
	D : \BZ^p \to \BZ^p; \quad (x_1,x_2,\ldots,x_p) \mapsto (|x_1-x_2|, |x_2-x_3|, \ldots, |x_n - x_1|),
	\]
	 satisfies $P|\frac{1}{3}p(2^{(p-1)/2}-1)$.
	\item (i) $\varepsilon_p \equiv 1 \bmod 2\CO_K$ or (ii) $3|h_p$.
\end{enumerate}
\end{Cor}

This strengthens the main result of \cite{Bre19}, in which only the implication $4(i)\Rightarrow (3)$ was shown. The equivalence $(1) \Leftrightarrow (3)$ is shown in \cite{BLM07}.

%%%%%%%%%%%%%%5
\section{Orders of fundamental units}

We first record the following result, see for example \cite[\S3]{Bro74}.
%\cite[pp. 98--99]{Wei84}.

\begin{Lem}\label{Norm-1}
If $p \equiv 1 \bmod 4$ is prime, then $N_{K/\BQ}(\varepsilon_p)=-1$ and $h_p$ is odd.
\end{Lem}

%\begin{Proof}
%%Sketch due to Keith Conrad on MathOverflow: Let $x + y\sqrt{p} > 1$ be the least unit greater than 1 with norm 1. Then $x$ and $y$ have to be positive integers (any unit in $\BZ[\sqrt{d}]$ which is greater than 1 must have positive coeff. wrt the basis $\{1,\sqrt{d}\}$). Show $x$ is odd, $y$ is even and from the equation $py^2 = (x+1)(x-1)$ get an equation $m^2 - pn^2 = 1$ where $1 < m < x.$ This can't go on forever, so some unit in $\BZ[\sqrt{p}]$ must have norm $-1$.
%
%\end{Proof}

%We collect here a few known results on the multiplicative order of the fundamental units $\varepsilon_p$ modulo $q$. 

The following result is due to Ishikawa and Kitaoka.

%When $q=2$, there are three possible residue classes for $\varepsilon_p \bmod 2 \in \BF_4^*$, and in the absence of any further structural results, we conjecture that each occurs equally often, hence $\varepsilon_p \equiv 1 \bmod 2$ for one third of all primes $p \equiv 5 \bmod 8$, see \cite[\S4]{Bre19} for numerical results supporting this. If one relaxes the condition $p$ prime to merely $p$ square-free, then results in this direction have been proved by Stevenhagen \cite{Stevenhagen96}.

\begin{Prop}\label{IK1}
Let $p \equiv 1 \bmod 4$ be a prime and suppose $q\nmid 2p$ is an inert prime in $K/\BQ$. Then 
\begin{enumerate}
\item $(q-1)/2$ divides $\ind(\varepsilon_p \bmod q)$
\item $\displaystyle \ord(\varepsilon_p \bmod q) \equiv 
\left\{\begin{array}{ll} 4 \bmod 8 & \text{if $q \equiv 1 \bmod 4$} \\
0 \bmod 8 & \text{if $q \equiv 3 \bmod 4$.} \end{array}\right.$
\end{enumerate}
In particular, $\ord(\varepsilon_p \bmod q) = \ord(-\varepsilon_p \bmod q)$.
\end{Prop}

\begin{Proof}
Since $N_{K/\BQ}(\varepsilon_p)=-1$ by Lemma \ref{Norm-1}, 
(1) follows from \cite[Theorem 1.1]{IK98}, and (2) follows from \cite[Corollary 1.4]{IK98}.
The final claim follows from the fact that $\ord(\varepsilon_p \bmod q)$ is divisible by~4.
\end{Proof}

Proposition \ref{IK1} and Theorem \ref{Thm:Main} imply that, if $q>3$, then $\ind(\alpha)$ is divisible by $(q-1)/2$. 

\newcolumntype{L}{>{\raggedright\arraybackslash}p{.4\textwidth}}

\begin{table}%
\renewcommand{\arraystretch}{1.5}
{\small
\begin{tabular}{|l|l|l|l|L|}
\hline
$q$ & $i(\varepsilon_p \bmod q)$ & Freq. pred. & Freq. obs. & Consequences of Theorem \ref{Thm:Main} \\
\hline\hline
2 & 1 & $\frac{2}{3}$ & 0.67497 & $3|\ind(\alpha)$ iff $3|h_p$  \\
\cline{2-5}
  & 3 & $\frac{1}{3}$ & 0.32503 & $3|\ind(\alpha)$ \\
\hline\hline
3 & 1 & 1 & 1.0 & $2\nmid\ind(\alpha)$ \\
\hline\hline
5 & 2 & $\frac{2}{3}$ & 0.67359 &  $2\|\ind(\alpha)$, and $3|\ind(\alpha)$ iff $3|h_p$ \\
\cline{2-5}
& 6 & $\frac{1}{3}$ &  0.32641 & $2\|\ind(\alpha)$ and $3|\ind(\alpha)$  \\

\hline\hline
7 & 3 & 1 & 1.0 & $2\nmid\ind(\alpha)$, and $3|\ind(\alpha)$ \\
\hline\hline
11 & 5 & $\frac{2}{3}$ & 0.67325 & $2\nmid\ind(\alpha)$, $5|\ind(\alpha)$, and $3|\ind(\alpha)$ iff $3|h_p$\\
\cline{2-5}
& 15 & $\frac{1}{3}$ & 0.32675 & $2\nmid\ind(\alpha)$ and $15|\ind(\alpha)$\\
\hline\hline
13 & 6 & $\frac{6}{7}$ & 0.85795 & $2\|\ind(\alpha)$, $3|\ind(\alpha)$, and $7|\ind(\alpha)$ iff $7|h_p$\\
\cline{2-5}
& 42 & $\frac{1}{7}$ & 0.14205 & $2\|\ind(\alpha)$ and $21|\ind(\alpha)$\\
\hline\hline
17 & 8 & $\frac{2}{3}$ & 0.67236 & $2^3\|\ind(\alpha)$, $3\|\ind(\alpha)$ iff $3\|h_p$, and $9|\ind(\alpha)$ iff  $9|h_p$ \\
\cline{2-5}
& 24 & $\frac{2}{9}$ & 0.21849 & $2^3\|\ind(\alpha)$, $3\|\ind(\alpha)$ iff $3\nmid h_p$, and $9|\ind(\alpha)$ iff  $3|h_p$\\
\cline{2-5}
& 72 & $\frac{1}{9}$ & 0.10914 & $2^3\|\ind(\alpha)$ and $9|\ind(\alpha)$ \\
\hline\hline
19 & 9 & $\frac{4}{5}$ & 0.80082 & $2\nmid\ind(\alpha)$, $9|\ind(\alpha)$ and $5|\ind(\alpha)$ iff $5|h_p$\\
\cline{2-5}
& 45 & $\frac{1}{5}$ & 0.19918 & $2\nmid\ind(\alpha)$ and $45|\ind(\alpha)$\\
\hline
\end{tabular}
}
\caption{Possible values of $\ind(\varepsilon_p \bmod q)$ for small $q$ and various $p$. The third and fourth columns list the predicted and observed frequency, respectively, of each given value of $\ind(\varepsilon_p \bmod q)$ for primes $p \equiv 1 \bmod 4, \; p < 10^8$.
}
\label{Table1}
\end{table}

The possible values of the indices $\ind(\varepsilon_p \bmod q)$ for $q\leq 19$ inert in $\BQ(\sqrt{p})/\BQ$ and $p \equiv 1 \bmod 4$ are listed in Table \ref{Table1}.
%, as are the consequences for $\ind(\alpha)$ and $\ind(\beta)$. 
%
These values are computed as follows.

Since $p\equiv 1 \bmod 4$ we have 
\[
\CO_K = \BZ\left[\frac{1+\sqrt{p}}{2}\right] \cong \frac{\BZ[X]}{\langle X^2 - X + \frac{1-p}{4}\rangle},
\]
and under this isomorphism,
\[
\varepsilon_p = \frac{x+y\sqrt{p}}{2} = \frac{1}{2}(x-y) + \left(\frac{1+\sqrt{p}}{2}\right)y \mapsto \frac{1}{2}(x-y) + Xy,
\]
where $x$ and $y$ satisfy the Pellian equation 
\begin{equation}\label{Pell}
x^2-py^2 = -4,
\end{equation}
since $N_{K/\BQ}(\varepsilon_p) = -1$ by Lemma \ref{Norm-1}.

Next, we consider the finite fields
\[
\CO_K/q\CO_K \cong \frac{\BF_q[X]}{\langle  X^2 - X + \frac{1-p}{4}\rangle},
\]
one for each residue class $p \bmod q$ such that $q$ is inert in $K/\BQ$. For odd $q$, we let $p$ range through the quadratic non-residues mod $q$, by quadratic reciprocity, and when $q=2$ we set $p=5$.
%
%
%When $q=2$ we use $p=5$, and when $q$ is odd we use 
%
%one for each quadratic non-residue class $p \bmod q$ when $q$ is odd (so $q$ is inert in $K/\BQ$, by quadratic reciprocity).

For each $a + bX$ in such a field, we check whether $a+bX = \frac{1}{2}(x-y) + Xy$ holds with $x,y$ satisfying (\ref{Pell}). If so, we compute its multiplicative index and we have found a candidate value for $\ind(\varepsilon_p \bmod q)$. The proportion of candidate residue classes for each multiplicative index gives a na\"ive prediction for the density of primes $p$ for which $\ind(\varepsilon_p \bmod q)$ equals that index. These predictions, together with the observed density for primes $p \equiv 1 \bmod 4, \; p < 10^8$ are shown in Table~\ref{Table1}. (Restricting to primes $p \equiv 5 \bmod 8$ produces similar results).

Lastly, the consequences for $\ind(\alpha)$ from Theorem \ref{Thm:Main} are also listed, using the fact that $h_p$ is odd.

\bigskip

We illustrate this with the example $q=5$. We consider the finite fields 
$F_2 = \BF_5[X]/\langle X^2-X-4 \rangle$ and $F_3 = \BF_5[X]/\langle X^2-X-3 \rangle$, corresponding to the residue classes $p \equiv 2$ and $3 \bmod 5$, respectively. In $F_2$ there are 6 elements $aX+b$ satisfying $(2a+b)^2-2b^2 \equiv -4 \bmod q$, which thus might represent $\varepsilon_p \bmod 5\CO_K$. Four of them, $2X, 3X, 2X+3$ and $3X+2$, have multiplicative index 2 while the elements $2$ and $3$ in $F_2$ have multiplicative index $6$. The situation is similar in $F_3$. Thus we predict that $\ind(\varepsilon_p \bmod 5)$ equals $2$ with probability $2/3$ and equals $6$ with probability $1/3$.

Theorem \ref{Thm:Main} says that $\gcd\big(\ind(\alpha),24\big) = \ind(\varepsilon_p^{h_p} \bmod 5)$.
It follows that $2\|\ind(\alpha)$, since $h_p$ is odd. Furthermore, $3|\ind(\alpha)$ if and only if $3|h_p$ or $\ind(\varepsilon_p \bmod 5) = 6$.

\bigskip

It would be interesting to prove that the predicted densities in Table~\ref{Table1} are indeed correct, but nothing seems to be known rigorously. Even the question of whether there are infinitely many primes $p \equiv 5 \bmod 8$ for which $\ind(\varepsilon_p \bmod 2) = 3$ is still open, although there are some known results if we relax the condition $p$ prime to $p$ squarefree, see \cite{Stevenhagen96}. 
%Some numerical evidence that this density should be $1/3$ is found in \cite[\S4]{Bre19}.

On the other hand, in the related situation in which $p$ is fixed and $q$ varies, more is known. In particular, densities of $q$ for which $\ind(\varepsilon_p \bmod q)$ equals a given value are obtained in \cite{CKY00,Kat03} under the assumption of the generalized Riemann Hypothesis.

%%%%%%%%%%%%%%%%%
\section{Proof of the main result}

From now on, we fix a prime number $p \equiv 1 \bmod 4$. 

Let $\zeta_p=\exp(2\pi i/p)\in\BC$ be a primitive $p^{\mathrm{th}}$ root of unity, $L=\BQ(\zeta_p)$ the corresponding cyclotomic number field and $L^+=\BQ(\zeta_p+\zeta_p^{-1})$ its maximal totally real subfield. 
Now $K=\BQ(\sqrt{p})$ is the unique quadratic subfield of $L^+$.

We denote by $\CO_L$ and $\CO_{L^+}$ the rings of integers of $L$ and $L^+$, respectively. 
 Our elements $\alpha\in\BF_{q^n}$ are the reductions of the unit $\zeta_p+\zeta_p^{-1}\in\CO_{L^+}^*$ modulo primes lying above $q$.

The following result determines the norms $N_{L/K}(\zeta_p+1)$ and $N_{L^+/K}(\zeta_p+\zeta_p^{-1})$  in~$K$.  It is one of many consequences of Dirichlet's Class Number Formula for~$K$; we prove it here for lack of a suitable reference.

\begin{Prop}
If $p\equiv 5 \bmod 8$ we have
\begin{enumerate}
	\item $\displaystyle N_{L/K}(\zeta_p+1) 
	= \prod_{k=1}^{(p-1)/2} \big(\zeta_p^{k^2}+1\big) 
	= \varepsilon_p^{2h_p}.$
		
	\item $\displaystyle N_{L^+/K}(\zeta_p+\zeta_p^{-1}) 
	= \prod_{0<r\leq\frac{p-1}{2}, \left(\frac{r}{p}\right)=1} \big(\zeta_p^{r}+\zeta_p^{-r}\big) 
	= (-1)^m\varepsilon_p^{h_p},$
	
	where
	\begin{align*}
	m &= \#\{r \;|\; \frac{p+3}{4} \leq r \leq \frac{p-1}{2}, \; \left(\frac{r}{p}\right)=1\} \\
	  &= \frac{1}{4}\left[ \frac{p-1}{2} - h(-p)\right]
	\end{align*}
	and $h(-p)$ is the class number of the imaginary quadratic field $\BQ(\sqrt{-p})$.
\end{enumerate} 
If $p \equiv 1 \bmod 8$ then $N_{L^+/K}(\zeta_p+\zeta_p^{-1})=(-1)^m$, with $m$ as above, and $N_{L/K}(\zeta_p+1)=1$.
\end{Prop}

\begin{Proof}
First note that
\begin{align*}
N_{L/K}(\zeta_p + 1) &= \prod_{\sigma\in\Gal(L/K)}\big(\sigma(\zeta)+1\big) \\
& = \prod_{0<r<p, \left(\frac{r}{p}\right)=1}\big(\zeta^r + 1\big) = \prod_{k=1}^{(p-1)/2} \big(\zeta_p^{k^2}+1\big).
\end{align*}
For an integer $a$ not divisible by $p$, we have the following formula, see
%which follows from the Class Number Formula, see for example 
\cite[Thm 1.3]{Sun2019}.
\begin{equation}
\label{eq:Sun}
\prod_{k=1}^{(p-1)/2} \big(1 - \zeta_p^{ak^2}\big) = \sqrt{p}\varepsilon_p^{-\left(\frac{a}{p}\right)h_p}
\end{equation}
Applying this with $a=1$ and $a=2$, we obtain
\begin{align*}
N_{L/K}(\zeta_p + 1) & = \prod_{k=1}^{(p-1)/2} \big(1 + \zeta_p^{k^2}\big) 
 = \frac{\prod_{k=1}^{(p-1)/2} \big(1 - \zeta_p^{2k^2}\big)}{\prod_{k=1}^{(p-1)/2} \big(1 - \zeta_p^{k^2}\big)} \\
& = \left\{\begin{array}{ll}
\varepsilon_p^{2h_p} & \text{if $p\equiv 5 \bmod 8$} \\
1 & \text{if $p \equiv 1 \bmod 8$.}
\end{array}\right.
\end{align*}
This proves (1). 

Next, 
\[ 
N_{L^+/K}(\zeta_p + \zeta_p^{-1}) = \prod_{\sigma\in\Gal(L^+/K)}\big(\sigma(\zeta_p) + \sigma(\zeta_p^{-1})\big) = \prod_r (\zeta_p^r + \zeta_p^{-r}),
\]
 where $r$ runs over exactly one non-zero quadratic residue $r \bmod p$ in each class modulo $\pm 1$, e.g. those quadratic residues $0<r\leq (p-1)/2$.

Furthermore,
\[
N_{L/K}(\zeta_p+\zeta_p^{-1}) = N_{L/K}\big(\zeta_p^{-1}(\zeta_p^2 + 1)\big)
= N_{L/K}(\zeta_p^2 + 1) = N_{L/K}(\zeta_p + 1),
\]
since $\zeta_p$ and $\zeta_p^2$ are conjugate and 
\[
N_{L/K}(\zeta_p)=\prod_{k=1}^{(p-1)/2}\zeta_p^{k^2} = \zeta_p^{p(p-1)(p+1)/24} = 1.
\]
Since $N_{L/K}(\zeta_p + \zeta_p^{-1}) = \big(N_{L^+/K}(\zeta_p + \zeta_p^{-1})\big)^2$, we have
obtained (2) up to a sign, which we determine next.

We have $\varepsilon_p^{h_p} > 0$, whereas the number of negative factors in $\prod(\zeta_p^r+\zeta_p^{-r}) = \prod 2\cos(2\pi i r/p)$ equals
\begin{align*}
m & = \#\{r \;|\; \frac{p}{4}<r<\frac{p}{2}, \; \left(\frac{r}{p}\right)=1\} \\
& = \frac{p-1}{4} - \#\{ r \;|\; 0<r<\frac{p}{4}, \; \left(\frac{r}{p}\right)=1   \} \\
& = \frac{p-1}{4} - \frac{1}{2}\sum_{r=1}^{(p-1)/4}\left[1+\left(\frac{r}{p}\right)\right] \\
& = \frac{p-1}{4} - \frac{1}{2}\left[ \frac{p-1}{4} + \sum_{r=1}^{(p-1)/4}\left(\frac{r}{p}\right) \right] \\
& = \frac{1}{4}\left[ \frac{p-1}{2} - h(-p) \right],
\end{align*}
by Dirichlet's class number formula for $h(-p)$, \cite[\S106]{Dir99}.
\end{Proof}

\begin{Rem}
We note that the formula in (2) above is similar to, but simpler than, the one discovered by P.~Chowla in \cite{Cho68}.
\end{Rem}

%\begin{Rem}
%If $p\equiv 3 \bmod 4$ then of course the unique quadratic subfield of $L$ is the imaginary quadratic $\BQ(\sqrt{-p})$, whose group of units is $\{\pm 1\}$ unless $p=3$.
%\end{Rem}

Now suppose $p \equiv 5 \bmod 8$ and $(\BZ/p\BZ)^* = \langle -1, q\rangle$. Then $q$ is inert in $L^+/\BQ$ and so the following diagram commutes, where the horizontal arrows are reduction modulo $q$ and the vertical arrows are norms.

\begin{equation}\label{eq:norms2}
\xymatrix@C=5pt{
\zeta_p+\zeta_p^{-1}\ar@{|->}[d] & \in & \CO_{L^+}^* \ar[rrr]\ar[d]^{N_{L^+/K}} &&& \big(\CO_{L^+}/q\CO_{L^+}\big)^* \ar@{->>}[d]^N \\
\pm\varepsilon_p^{h_p} & \in & \CO_K^* \ar[rrr] &&& \big(\CO_K/q\CO_K)^*
}
\end{equation}
Here, $N(\alpha) = \pm\varepsilon_p^{h_p} \bmod q\CO_K$, where the sign is irrelevant for the multiplicative index by Proposition \ref{IK1}.
% when $q>2$, the case $q=2$ being clear.

Theorem \ref{Thm:Main} now follows from the following lemma, where for an element $g$ of a finite group $G$, we denote by $\ord_G(g)$ its order and by $\ind_G(g) = \#G/\ord_G(g)$ its index.

\begin{Lem}
Let $f:G\onto H$ be an epimorphism of finite cyclic groups and $g\in G$. Then $\ind_H(f(g))=\gcd(\ind_G(g),|H|)$
\end{Lem}

\begin{Proof}
For every divisor $d$ of $|H|$, denote by $H_d < H$ and $G_d < G$ the unique subgroup of index $d$.
Let $\ell$ be a prime number dividing $|H|$ and let $n=v_{\ell}(|H|)$ be the $\ell$-adic valuation of $|H|$. Then for every $0\leq i \leq n$, the map $f$ restricts to an epimorphism $f : G_{\ell^i} \onto H_{\ell^i}$. Now 
\begin{align*}
v_{\ell}(\ind_G(g)) = m & \quad \Longleftrightarrow \quad g\in G_{\ell^m}\smallsetminus G_{\ell^{m+1}} \\ 
& \quad \Longleftrightarrow \quad f(g)\in H_{\ell^m}\smallsetminus H_{\ell^{m+1}} \\
& \quad \Longleftrightarrow \quad v_{\ell}(\ind_H(f(g))) = m.
\end{align*}
The result follows.
\end{Proof}

%%%%%%%%%%%%
\section{Some heuristics}

How often does a given prime divide $\ord(\alpha)$?

Suppose $d|q^n-1$. Then a randomly chosen element $\beta\in\BF_{q^n}^*$ satisfies $d|\ind(\beta)$ with probability $1/d$, since $\BF_{q^n}^*$ has a unique subgroup of index~$d$. 

In the case $q=2$ and $p\equiv 5 \bmod 8$, a na\"ive heuristic (e.g. \cite[\S4]{Bre19}) suggests that 
$\varepsilon_p \equiv 1 \bmod 2$ occurs with probability $1/3$, whereas the Cohen-Lenstra heuristics \cite[\S9.II]{CohenLenstra} predict that $3|h_p$ with probability $1 - \prod_{k\geq 2}(1-3^{-k})\approx 0.159811$. Assuming that these conditions  are independent, we thus expect the index $\ind(\zeta + \zeta^{-1})$ to be divisible by $3$ for about $43.9874\%$ of primes $p\equiv 5 \bmod 8$ for which $2$ is a primitive root. 

This suggests that the Gauss period $\alpha = \zeta + \zeta^{-1} \in\BF_{2^n}^*$ is at least $10\%$ less likely to be a primitive root than a randomly chosen element, due to the potential $3$-divisibility of the class number $h_p$.

\medskip

Lastly, we consider the case where $p=2r+1$ and $r$ is also prime, in which case $r$ is called a Sophie Germain prime. Since $[\BQ(\zeta_p + \zeta_p^{-1}) : \BQ] = r$ is prime, there are no intermediate fields $K$ for which a phenomenon like Theorem \ref{Thm:Main} might occur. In this case, a conjecture of Gao and Vanstone \cite{GV95} states that the Gauss period $\alpha \in \BF_{2^r}^*$ is always a primitive root. The conjecture is verified in \cite{GV95} for $r < 593$.

We give some heuristic arguments supporting this conjecture.

Assuming the Gauss period $\alpha$ behaves like a random element of $\BF_{2^r}^*$, any prime divisor $\ell$ of $2^r-1$ will divide $\ind(\alpha)$ with probability $1/\ell$. What is the probability that $\ell$ divides $2^r-1$? Na\"ively, we expect $1/\ell$ also. Less Na\"ively, we may argue as follows (see e.g. \cite{Wag83}).

Every prime divisor $\ell | 2^r-1$ (where $r$ is prime) must be of the form $\ell = 2kr+1$, where $k \equiv 0$ or $-r \bmod 4$. The proportion of primes of this form is $\frac{1}{2}\frac{1}{\phi(2r)} \approx \frac{1}{2r}$. By a heuristic argument from \cite{SK67}, each such prime has probability $1/k \approx 2r/\ell$ of dividing $2^r-1$. Combining these, we again find that a prime $\ell > 2r$ divides $2^r-1$ with probability $1/\ell$.

%Define $G(x) := \int_{x}^\infty \frac{dt}{t^2\log t}$, so $\sum_{\text{$\ell>x$ prime}} 1/\ell^2 \approx G(x)$.

The expected number of counter-examples to the conjecture of Gao and Vanstone is thus less than
\begin{align*}
&\sum_{\scriptsize \begin{array}{c}r\geq 593 \\ \text{$r$ Sophie Germain prime}\end{array}} 
\sum_{\scriptsize \begin{array}{c}\ell>2r \\ \text{$\ell$ prime}\end{array}}\frac{1}{\ell^2} 
 \approx \sum_{r=593}^\infty \frac{2C}{\log^2 r} \sum_{l = 2r+1}^\infty \frac{1}{l^2\log l} \\
& \approx\int_{593}^\infty \frac{2C}{\log^2 r} \left(\int_{2r+1}^\infty \frac{1}{l^2\log l} \,dl\right) dr
\approx 0.007.
\end{align*}

Here we have used the heuristic that $r$ is a Sophie Germain prime with probability $2C/\log^2 r$, where $C\approx 0.66$ is the Hardy-Littlewood twin prime constant.

%\begin{align*}
%\; &\sum_{\scriptsize \begin{array}{c}r\geq 593 \\ \text{$r$ prime} \\ \text{$2r+1$ prime}\end{array}} \left(1-\prod_{\scriptsize \begin{array}{c}\ell>2r \\ \text{$\ell$ prime}\end{array}}\left(1-\frac{1}{\ell^2}\right)\right) \\
%& \approx \sum_{\scriptsize \begin{array}r\geq 593\\ \text{$r$ prime}\end{array}} \frac{1.32}{\log(2r+1)}
 %\sum_{\scriptsize \begin{array}{c}\ell>2r \\ \text{$\ell$ prime}\end{array}}\frac{1}{\ell^2} \\
%& < 1.32\sum_{n=195}^\infty \frac{1}{p_n^2}\left[\sum_{k=195}^n\frac{1}{\log p_k}\right]  \\
%& < 1.32\sum_{n=195}^\infty \frac{n}{p_n^2\log n} 
 %< 1.32\sum_{n=195}^\infty \frac{1}{n\log^3 n} \\
%& < 1.32\int_{194}^\infty \frac{dv}{v\log^3 v} < 0.024.
%\end{align*}

%
%Here we have used the fact \cite{Ros39} that the $n$th prime $p_n > n\log n$, as well as the heuristic that a prime $r$ is a Sophie Germain prime (i.e. $2r+1$ is also prime) with probability  $\frac{2C}{\log r}$, where $C\approx 0.66$ is the Hardy-Littlewood twin prime constant.

\subsection*{Acknowledgements.}
The author would like to thank Igor Shparlinski for helpful comments.


\begin{thebibliography}{99}

\bibitem[ASV10]{ASV10}
O. Ahmadi, I. E. Shparlinski and J. F. Voloch,
Multiplicative order of Gauss periods,
{\em Internat. J. Number Theory} {\bf 6} (2010) no. 4, 877--882.

%\bibitem[And76]{And76}
%G. E. Andrews,
%The Theory of Partitions,
%Addison-Wesley, New York, 1976.

%\bibitem[Ar06]{Arnold06} 
%V. I. Arnold, 
%Complexity of finite sequences of zeroes and ones and geometry of finite spaces of functions, 
%{\em Funct. Anal. Other Math.} {\bf 1} (2006), 1--15. 

%\bibitem[Av13]{Av13}
%C. Avart, 
%A result about cycles in Ducci sequences, 
%{\em Fibonacci Quart.} {\bf 51} (2013), no. 2, 137--141.

%\bibitem[BKP05]{BKP} 
%A. Behn, C. Kribs-Zaleta, V. Ponomarenko, 
%The convergence of difference boxes, 
%{\em Amer. Math. Monthly} {\bf 112 (5)} (2005) 426--439.

%\bibitem[Ber75]{Berlekamp75} 
%E. R. Berlekamp, 
%The design of slowly shrinking labelled squares, 
%{\em Math. Comp.} {\bf 29} (1975) 25--27.




\bibitem[BLM07]{BLM07} F. Breuer, E. L\"otter and A.B. van der Merwe,  
Ducci sequences and cyclotomic polynomials, 
{\em Finite Fields Appl.} {\bf 13} (2007), 293--304.


%\bibitem[Bre10]{Bre10}  
%F. Breuer, 
%Ducci sequences and cyclotomic fields, 
%{\em J. Difference Equ. Appl.} {\bf 16} (2010), no. 7, 847--862.

\bibitem[Bre19]{Bre19}
F. Breuer,
Periods of Ducci sequences and odd solutions to a Pellian equation,
{\em Bull. Aust. Math. Soc.}, {\bf 100} (2019), 201--205
doi:10.1017/S0004972719000212

\bibitem[BS19]{BS19}
F. Breuer and I. E. Shparlinski,
Lower bounds for periods of Ducci sequences,
{\em Bull. Aust. Math. Soc.}, To appear.

%\bibitem[BZ07]{BZ07} 
%G. Brockman, R. J. Zerr, 
%Asymptotic behaviour of certain Ducci sequences, 
%{\em Fibonacci Quart.} {\bf 45} (2007), no. 2, 155--163.

\bibitem[Bro74]{Bro74}
E. Brown, 
Class numbers of real quadratic number fields.
{\em Trans. Amer. Math. Soc.} {\bf 190} (1974), 99--107.

%\bibitem[BM08]{BM08}
%R. Brown, J. L. Merzel, 
%The number of Ducci sequences with given period, 
%{\em Fibonacci Quart.} {\bf 45} (2007), no. 2, 115--121. 

%\bibitem[BFJ78]{BFJ78} 
%M.~Burmester, R.~Forcade and E.~Jacobs, 
%Circles of numbers, {\em Glasgow Math. J.} {\bf 19} (1978), 115--119.

\bibitem[CKY00]{CKY00}
Y.-M. Chen, Y. Kitaoka and J. Yu,
Distribution of units of real quadratic number fields,
{\em Nagoya Math. J.} {\bf 158} (2000), 167--184.

\bibitem[Cho68]{Cho68}
P. Chowla, 
On the class-number of real quadratic fields,
{\em J. Reine Angew. Math.} {\bf 230} (1968), 51--60.

%\bibitem[CM37]{CM37} 
%C.~Ciamberlini and A.~Marengoni, 
%Su una interessante curiosit\`a numerica,
%{\em Periodiche di Matematiche} {\bf 17} (1937), 25--30.

%\bibitem[CST05]{CST05}
%N. J. Calkin, J. G. Stevens, D. M. Thomas, 
%A characterization for the length of cycles of the n-number Ducci game,
%{\em Fibonacci Quart.} {\bf 43} (2005), no. 1, 53--59. 

%\bibitem[Cla18]{Clausing} 
%A.~Clausing, 
%Ducci Matrices, 
%{\em Amer. Math. Monthly}, {\bf 125} (2018), no. 10, 901--921.

\bibitem[CL84]{CohenLenstra}
H. Cohen and H. W. Lenstra,
Heuristics on class groups of number fields,
in: {\em Number Theory Noordwijkerhout 1983}, 33--62, Springer-Verlag, 1984.


\bibitem[Dir99]{Dir99}
P.G.L. Dirichlet,
Lectures on Number Theory,
{\em History of Mathematics Sources} {\bf 16}, Amer. Math. Soc. \& London Math. Soc., 1999.


%\bibitem[Ehr90]{Ehr90}
%A. Ehrlich, 
%Periods of Ducci's N-number game of differences, 
%{\em Fibonacci Quart.} {\bf 28} (1990), no. 4, 302--305.

%\bibitem[Fre48]{Freedman48} 
%B. Freedman, 
%The four number game, 
%{\em Scripta Math.} {\bf 14} (1948) 35--47, reprinted in arXiv:1109.0051v1.


%\bibitem[GS95]{GS95}
%H. Glaser and G. Sch\"offl,
%Ducci-sequences and Pascal's triangle,
%{\em Fibonacci Quart.} {\bf 33} (1995) No. 4, 313--324.

\bibitem[GS98]{GS98}
J. von zur Gathen and I. E. Shparlinski,
Orders of Gauss periods in finite fields,
{\em . Appl. Algebra Engrg. Comm. Comput.} {\bf 9} (1998), no. 1, 15--24.

\bibitem[GV95]{GV95}
S. Gao and S. A. Vanstone,
On orders of optimal basis generators,
{\em Math. Comp.} {\bf 64} (1995), no. 211, 1227--1233. 

%\bibitem[Hag64]{Hag64}
%P. Hagis,
%On a Class of Partitions with Distinct Summands,
%{\em Trans. Amer. Math. Soc.} {\bf 112} (1964), No. 3, 401--415.

\bibitem[IK98]{IK98}
M. Ishikawa and Y. Kitaoka,
On the distribution of units modulo prime ideals in real quadratic fields,
{\em J. reine angew. Math.} {\bf 494} (1998), 65--72.

\bibitem[Kat03]{Kat03}
N. Kataoka,
The distribution of prime ideals in a real quadratic field with units having a given index in the residue class field,
{\em J. Number Theory} {\bf 101} (2003), no. 2, 349--375.

%\bibitem[Lud81]{Lud81}
%A. L. Ludington,
%Cycles of differences of integers,
%{\em J. Number Theory}, {\bf 13} (1981), 255--261.

%\bibitem[MST06]{MST06}
%M. Misiurewicz, J. G. Stevens, D. M. Thomas, 
%Iterations of linear maps over finite fields, 
%{\em Linear Algebra Appl.} {\bf 413} (2006), no. 1, 218--234. 

%\bibitem[OEIS]{OEIS} On-line Encyclopedia of Integer Sequences, 
%entry $\#$A038553. https://oeis.org/A038553.

\bibitem[Pop12]{Pop12}
R. Popovych, 
Elements of high order in finite fields of the form $\BF_q[x]/\Phi_r(x)$, 
{\em Finite Fields Appl.} {\bf 18} (2012), No. 4, 700--710.

\bibitem[Pop14]{Pop14}
R. Popovych,
Sharpening of the explicit lower bounds for the order of elements in finite field extensions based on cyclotomic polynomials. 
{\em Ukrainian Math. J.} {\bf 66} (2014), no. 6, 916--927.

%\bibitem[Ros39]{Ros39}
%B. Rosser,
%The $n$th prime is greater than $n\log(n)$.
%{\em Proc Lond. Math. Soc. (2)}, {\bf 45} (1939), 21--44.

\bibitem[SK67]{SK67}
D. Shanks and S. Kravitz,
On the distribution of Mersenne divisors.
{\em Math. Comp.} {\bf 21} (1967), 97--101.


%\bibitem[Sim93]{Sim93} 
%G. J.~Simmons, 
%The structure of the differentiation digraphs of binary sequences, 
%{\em Ars Combin.} {\bf 35} (1993), A, 71--88.

%\bibitem[SB18]{SB18}
%S. Solak, M. Bah\c{s}i, 
%Some properties of circulant matrices with Ducci sequences, 
%{\em Linear Algebra Appl.} {\bf 542} (2018), 557--568.

\bibitem[Ste96]{Stevenhagen96} 
P. Stevenhagen, 
On a problem of Eisenstein, 
{\em Acta Arith.} {\bf 74} (1996), no. 3, 259--268.

\bibitem[Sun19]{Sun2019}
Z.-W. Sun,
Quadratic residues and related permutations and identities,
{\em Finite Fields Appl.} {\bf 59} (2019), 246--283.

\bibitem[Wag83]{Wag83}
S.S. Wagstaff, Jr.,
Divisors of Mersenne Numbers,
{\em Math. Comp.} {\bf 40} (1983), 385--397.

%\bibitem[Wei84]{Wei84}
 %A. Weil, 
%Number Theory: An Approach through History, 
%Birkh\"auser, Boston (1984).

%\bibitem[Zve79]{Zven79} 
%P.~Zvengrowski, 
%Iterated absolute differences, 
%{\em Math. Mag.} {\bf 52 (1)} (1979), 36--40.


\end{thebibliography}
\end{document}